\documentclass[sn-mathphys,Numbered]{sn-jnl}

\usepackage{graphicx}%
\usepackage{multirow}%
\usepackage{amsmath,amssymb,amsfonts}%
\usepackage{amsthm}%
\usepackage{mathrsfs}%
\usepackage[title]{appendix}%
\usepackage{xcolor}%
\usepackage{textcomp}%
\usepackage{manyfoot}%
\usepackage{tabularx}
\usepackage{multirow}
\usepackage{booktabs}%
\usepackage{algorithm}%
\usepackage{algorithmicx}%
\usepackage{algpseudocode}%
\usepackage{listings}%
\usepackage{tabularx}
\usepackage{fullwidth}

\theoremstyle{thmstyleone}%
\newtheorem{theorem}{Theorem}
\newtheorem{proposition}[theorem]{Proposition}%

\theoremstyle{thmstyletwo}%
\newtheorem{remark}{Remark}%
\newtheorem{lemma}{Lemma}
\newtheorem{corollary}{Corollary}

\theoremstyle{thmstylethree}%
\newtheorem{definition}{Definition}%

\usepackage{csquotes}
\usepackage{float}
\usepackage{mathtools}

\usepackage{tabularx}
\usepackage{multirow}
\usepackage{fullwidth}

\usepackage{csquotes}
\usepackage{float}
\usepackage{mathtools}

\usepackage{caption}
\begin{document}

\title[Article Title]{Statistical complexity as a criterion for the useful signal detection problem}


\author*[1]{\fnm{Leonid} \sur{Berlin}}\email{berlin.lm@phystech.edu}
\equalcont{These authors contributed equally to this work.}

\author[1]{\fnm{Andrey} \sur{Galyaev}}\email{galaev@ipu.ru}
\equalcont{These authors contributed equally to this work.}

\author[1]{\fnm{Pavel} \sur{Lysenko}}\email{pavellysen@ipu.ru}
\equalcont{These authors contributed equally to this work.}

\affil[1]{\orgdiv{Laboratory 38}, \orgname{Institute of Control Sciences of RAS}, \orgaddress{\city{Moscow}, \country{Russia}}}


\abstract{Three variants of the statistical complexity function, which is used as a criterion in the problem of detection of a useful signal in the signal-noise mixture, are considered. The probability distributions maximizing the considered variants of statistical complexity are obtained analytically and conclusions about the efficiency of using one or another variant for detection problem are made.
The comparison of considered information characteristics is shown and analytical results are illustrated on an example of synthesized signals. A method is proposed for selecting the threshold of the information criterion, which can be used in decision rule for useful signal detection in the signal-noise mixture. The choice of the threshold depends a priori on the analytically obtained maximum values. As a result, the complexity based on the total variation demonstrates the best ability of useful signal detection.}




\maketitle

\section{Introduction}\label{Intro}

The concept of information entropy was firstly introduced in Claude Shannon's article \cite{Shannon} in 1948. This work marked the beginning of a new field of science called information theory \cite{EntInf}. The development of information theory made possible an analytical and practical research in many applied fields of science and technology. Such terms as Gibbs and von Neumann entropies, Kullback-Leibler distance, Jensen-Shannon divergence, information divergences and some others were introduced and interpreted and later began to serve as criteria for various optimization problems of recognition \cite{recognition}, classification \cite{classification} and filtering.

By the end of the last century various information criteria, mainly Shannon information entropy, had began to be actively applied in the tasks of digital signal processing, in particular in the problem of detection of a useful signal in a noise environment \cite{RobustDetection}. The concept of spectral entropy \cite{EntUni}, associated with the Fourier spectrum of the considered signal, has appeared and proved to be especially relevant in the analysis of acoustic signals \cite{VAD}. In addition, the entropic approach has been successfully applied in the analysis of time series in the medical field, such as ECG or EEG \cite{ECG}. Later, a statistical complexity function was proposed as a development of the entropy concept \cite{ComplexIntro, Complexity, Complex_???}. However, the articles mostly do not provide an analytical study of the properties of these functions, which turns out to be especially important when solving the problem of hypothesis testing. 

It should be noted that there are several classical ways of solving the detection problem. The first of them is based on solving the problem of optimal filtering and requires knowledge of the properties of the signal: periodicity, bandwidth, etc. \cite{SignalAnalysis}. The second way is based on the Neumann-Pearson Lemma, solves the problem of hypothesis testing, and determines the fact of exceeding the optimal threshold at a given false alarm probability and requires estimation of statistical properties of sample distributions of noise and mixture of signal and noise \cite{Shiryaev}. The third way is equivalent to solving the changepoint detection problem when the unknown statistical characteristics of the signal distributions change. The anomaly detection problem \cite{Anomaly} has a similar formulation. All these methods demonstrate qualitative and reliable performance when the signal exceeds the noise, but for small signal-to-noise ratios often give the wrong answer.

The article is devoted to the problem of detection of useful signal in the signal-noise mixture and combines all three previously listed ways of solving the detection problem. We propose to use a variant of the Neumann-Pearson Lemma for the problem of hypothesis testing \cite{Shiryaev}, which is indeed valid when the error probability is close to one and depends on the total variation of the measure of two distributions of the null and alternative hypotheses. Based on the analytical expression of this error function, the criterion of the signal detection problem is formalized as one of the variants of the statistical complexity \cite{Sensors}, which takes into account the deterministic nature of the signal mixed in with the noise. The peculiarity of the statistical complexity is that it is multiplicative and consists of two multipliers, one of which is zero on deterministic sinusoidal signals of the same frequency (in physics these are objects of a given structure, such as crystals \cite{ComplexIntro}) and the other is zero on uniform distribution functions \cite{Complexity}, corresponding, for example, to white noise. Then the introduced criterion is compared with the already known two variants of statistical complexity based on Euclidean distance square and Jensen-Shannon divergence, their properties are established, and optimization as a function of many variables on a set of discrete distributions is performed. As a result, families of optimal distributions are identified and maxima of statistical complexity functions are calculated.

The article has the following structure. Section \ref{Intro} provides a literature review and highlights the current state of research on the topic of the article. Section \ref{hypothesis} is devoted to the connection of the considered information criteria  with the classical criterion of the signal detection problem. Section \ref{statoptimize} investigates the properties of the three types of statistical complexity. In \ref{Modelling} the analytical results of the previous section are supported by numerical simulations for synthesized signals. Section \ref{Conclusion} summarizes the results obtained in the paper and lists plans for the future.

\section{Neyman-Pearson Lemma and statistical complexity}\label{hypothesis}

The problem of signal detection $s(n)$ is traditionally reduced to the problem of hypothesis testing
\begin{equation*}
    \left\{
    \begin{array}{lll}
    \Gamma_0: x(n) = w(n),\\
    \Gamma_1: x(n) = s(n) + w(n),~n = 1, \dots, N.
    \end{array}
    \right.
\end{equation*}
Hypothesis $\Gamma_0$ corresponds to the decision of receiving only noise, and hypothesis $\Gamma_1$ -- of receiving a mixture of useful signal and noise, where the sequences $\{x(n)\},~n = 1, \dots, N$ are time series of the received data, $\{s(n)\}$ -- useful signal, $\{w(n)\}$ --additive white Gaussian noise, $N$ -- the length of the time series of data.

The random variables of the time series $(x(1),\dots,x(n),\dots,x(N))$ take values $(x_1,\dots,x_n,\dots,x_N) \in \mathbb{R}^N$.
 
In order to obtain an analytical expression for estimating the error probability in hypothesis testing, we can apply a variant of the Neyman-Pearson Lemma \cite{Shiryaev, DetectingDook}.

\begin{lemma} [Neyman-Pearson] \label{NeymanPearson}
Let there be an arbitrary, called a decision rule or test, measurable function of many variables \mbox{$(x_1, \dots, x_N) \in \mathbb{R}^N$} such that
\begin{equation*} 
d(x_1, ..., x_N)= \left\{{}
\begin{array}{l}
    1, ~\textrm {hypothesis}~\Gamma_0 ~\textrm {is true},\\\\
    0, ~\textrm {hypothesis}~\Gamma_1 ~\textrm {is true}, 
\end{array}
\right.
\end{equation*}
by which the following probabilities can be determined:

\begin{equation*}
    \begin{array}{cc}
    \alpha(d)=\textrm {Probability~(accept~} \Gamma_0 | \Gamma_1~\textrm {is true}),\\\\
    \beta(d)=\textrm {Probability~(accept~} \Gamma_1 | \Gamma_0~\textrm {is true}).
    \end{array}
\end{equation*}

Then the decision rule $d^*$ is optimal if

\begin{equation}
\alpha(d^*)+\beta(d^*)=\inf_d[\alpha(d)+\beta(d)]=\mathcal{E}r(N;\Gamma_0,\Gamma_1)~\text{-- error function},
\end{equation}
where the infinum is taken for all tests.
\end{lemma}
Here $\alpha(\cdot)$ is the probability of a false alarm, and $\beta(\cdot)$ is the probability of a useful signal missing.

The exact formula for the error function is as follows:
\begin{equation}\label{ErrTVR}
\displaystyle \mathcal{E}r(N;\Gamma_0,\Gamma_1)=1-\frac{1}{2}\lVert P_0^{(N)}-P_1^{(N)}\rVert = 1 - TV(P_0, P_1),
\end{equation}
where $P_0^{(N)}$ is the multivariate distribution function of the observation statistics by hypothesis $\Gamma_0$, $P_1^{(N)}$ is the multivariate distribution function of the observation statistics by hypothesis $\Gamma_1$, and $TV(P_0, P_1)$ is the total variation of the signed measure, \mbox{$\|Q\|=2\sup_A|Q(A)|$}. Thus, if the supports of measures $P_0$, $P_1$ do not overlap, then error-free distinguishing of hypotheses is possible. If the measures $P_0^{(N)}$ and $P_1^{(N)}$ are close, then $\lVert P_0^{(N)}-P_1^{(N)}\rVert \approx0$, leading to $\mathcal{E}r(N;\Gamma_0,\Gamma_1)\approx 1$.

For the problem of detecting a deterministic useful signal, for example, at a small signal-to-noise ratio, the case $\lVert P_0^{(N)}-P_1^{(N)}\rVert = 2 TV(P_0, P_1) \approx 0$ is of interest and the possibility of reasonable estimation of this value. Therefore, when the probability of the total error of distinguishing two hypotheses is close to one, it becomes possible to use the analytical expression $TV(P_0, P_1)$ to design a criterion in the problem of detecting a useful signal in a mixture. But first let us turn to already known criteria and establish their properties.

Most often, for the convenience of mathematical investigation, both of the useful signal and noise are modeled by Gaussian random processes with different parameters. In that case the problem of finding the moment of appearance of the signal $s(n)$ in the received sequence of samples is called the the problem of changepoint detection \cite{Shiryaev}.

Here and below, we consider discrete probability distributions $p = (p_1,\dots,~p_i,\dots,p_N)$, that by definition have the following properties:
\begin{equation}\label{discdefinit}
    \forall~ p_i \in [0, 1], \quad \sum_{i=1}^N p_i = 1.
\end{equation}

To formalize criterion that takes into account the deterministic component of the signal as well as the random one, let us explore the concepts of disequilibrium function $D$ and statistical complexity $C$ of the distribution. The simplest example of the disequilibrium function  is the square of Euclidean distance in the space of discrete probability distributions\cite{Complexity}.
\begin{definition}\label{Def1}
The disequilibrium $D_{SQ}$ has the meaning of the variance of a distribution relative to a uniform distribution
\begin{equation}\label{disequ_1}
    \displaystyle D_{SQ}(p)= \sum_{i=1}^N\left(p_i - \frac{1}{N} \right)^2=\sum_{i=1}^N p_i^2 - \frac{1}{N} .
\end{equation}
\end{definition}

\begin{definition}\label{Def2}
The statistical complexity, defined through the expression of disequilibrium by the Definition \ref{Def1}, is equal to
\begin{equation}\label{compSQ}
   C_{SQ}(p)= H(p)\cdot D_{SQ}(p),
\end{equation}
where
\begin{equation}\label{shannon_entropy}
    H(p) = \frac{1}{\log N}\left(-\sum_i^N p_i\log p_i \right)
\end{equation}
-- Shannon's normalized entropy \cite{Shannon}.
\end{definition}

In evaluating the sum \eqref{shannon_entropy}, it is assumed that $\displaystyle \frac{0}{\log 0} = 0$ by continuity, and this assumption holds for all subsequent equations.

It follows from the Definition \ref{Def1} that disequilibrium of the form (\ref{disequ_1}) and complexity of the form (\ref{compSQ}) are convenient to apply in estimation and comparison of signals having spectral distribution close to uniform. In general, instead of a uniform distribution $q_i=1/N$ at $i=1,...,N$, the formula (\ref{disequ_1}) may include an arbitrary discrete distribution. 

The formula \eqref{disequ_1} is proposed in \cite{Complexity} for computing the disequilibrium with respect to a uniform distribution, but most studies use the Jensen-Shannon divergence $JSD(p || q)$ \cite{Features} instead.

\begin{definition} \label{Def3}
The Jensen-Shannon disequilibrium equals
\begin{equation}\label{DJSD}
D_{JSD}(p)=JSD(p||q),
\end{equation}
where $q=(1/N, \dots, 1/N)$ is the uniform distribution.
\end{definition}
\begin{definition}\label{Def4}
Statistical complexity defined through the expression of disequilibrium from Definition \ref{Def3}, is expressed as
\begin{equation}\label{compJSD}
    C_{JSD}(p) = H(p)\cdot D_{JSD}(p).
\end{equation}
\end{definition}

\begin{remark}
It was noted above that $\displaystyle\sqrt{D_{SQ}}$ is a Euclidean metric on the space of discrete distributions. At the same time $\displaystyle\sqrt{D_{JSD}}$ is also a metric which is proportional to the Fisher metric.
\end{remark}

Since the error function of distinguishing between two hypotheses depends on the total variation $TV(p,q)$, which is obtained in the Neyman-Pearson Lemma \ref{NeymanPearson}, we introduce another notion of disequilibrium.

\begin{definition}\label{Def5}
The disequilibrium based on the total variation of signed measure is equal to
\begin{equation}\label{DTV}
D_{TV}(p)=TV^2(p,q),
\end{equation}
where $q=(1/N, \dots, 1/N)$.
\end{definition}
\begin{definition}\label{Def6}
The statistical complexity, defined through the disequilibrium expression according to the Definition \ref{Def5}, is equal to
\begin{equation}\label{compTV}
    C_{TV}(p) = H(p)\cdot D_{TV}(p).
\end{equation}
\end{definition}

The information divergence functions presented above, which define different variants of the disequilibrium function, can be unified by the general concept of \textit{f--divergence} \cite{fdivergence}:
\begin{equation}\label{generalfdiv}
    D_f(p||q) = \sum_{x \in \mathbb{R}^N} q(x) f\left(\frac{p(x)}{q(x)}\right).
\end{equation}

The choice of function $f$ gives rise to a whole family of different divergences:

\begin{itemize}
    \item The Kulback-Leibler divergence $D_{KL}(p, q)$ is obtained from \eqref{generalfdiv} by choosing $f(x) = x\log(x),~x > 0$.
    \item The Jensen-Shannon divergence is obtained from \eqref{generalfdiv} by choosing 
    \begin{equation}
        \displaystyle f(x) = x\log\frac{2x}{x+1} + \log \frac{2}{x+1}, ~x > 0.
    \end{equation}
     \item The total variation is obtained when $\displaystyle f(x) = \frac{1}{2}|1 - x|$:
    \begin{equation}\label{TVequ}
        \displaystyle TV(p, q)={\frac {1}{2}}\sum _{x \in \mathbb{R}^N }|p({x)-q(x)|};
    \end{equation}
 $TV(p,q)$ is also a metric on the space of probability distributions.
    The total variation is related to the Jensen-Shannon divergence by the following relation:
    \begin{equation}\label{TVvsJSD}
        JSD(p || q) \le TV(p,q).
    \end{equation}
\end{itemize}

It follows from the inequality (\ref{TVvsJSD}) that the total variation is the upper bound of Jensen-Shannon divergence.

Next, let us investigate the possibility of using each variant of statistical complexity as a criterion for indicating the appearance of a signal, but at first their properties must be established.

\section{Statistical complexity optimization}\label{statoptimize}
\subsection{Optimization of $C_{SQ}$}
Let us formulate the problem of maximizing the statistical complexity function on the set of discrete distributions $p=(p_1,...,p_N)$
\begin{equation}\label{comp}
   C_{SQ}(p)= \frac{1}{\log N}\left(-\sum_{i=1}^N p_i\log p_i \right)\cdot \left(\sum_{i=1}^N \left(p_i - \frac{1}{N} \right)^2\right)\longrightarrow \max_{p}
\end{equation}
with the condition
\begin{equation}\label{norm}
    \sum_{i=1}^N p_i = 1.
\end{equation}

An auxiliary result will be needed to formulate the Lemma about the maximum value of statistical complexity.

\begin{lemma}\label{lemma3}
Let $ 0< x\leq y \leq z\leq 1$, then $f(x,y,z)=x^y y^{-x}z^x x^{-z}y^z z^{-y}\geq 1$, with equality possible only when either $x=y$ or $y=z$.
\end{lemma}
\begin{proof}
Let us introduce a new function $g(x,y,z)=\ln f(x,y,z)$, 
$$
g(x,y,z)=y\ln x-x\ln y+x\ln z-z\ln x+z\ln y-y\ln z.
$$
Then it is required to prove that $g(x,y,z)\geq 0$ for $ 0 <x\leq y \leq z\leq 1$.

By the Kuhn-Tucker theorem, the solution of the conditional optimization problem of a function of three variables is either at the interior point of the constraint manifold or at its boundary.
The necessary conditions for the unconditional extremum of the function $g(x,y,z)$ take the following form
\begin{equation}\label{ENC}
\begin{array}{l}
\displaystyle \frac{\partial g}{\partial x}=\ln z-\ln y+\frac{y-z}{x}=0,\\
\displaystyle \frac{\partial g}{\partial y}=\ln x-\ln z+\frac{z-x}{y}=0,\\
\displaystyle \frac{\partial g}{\partial z}=\ln y-\ln x+\frac{x-y}{z}=0.
\end{array}
\end{equation}
Let us summarize all the equations of the last system:
\begin{equation*}
\frac{y-z}{x}+\frac{z-x}{y}+\frac{x-y}{z}=0,
\end{equation*}
which can be rewritten as
\begin{equation*}
\frac{(y-z)(x-y)(z-x)}{xyz}=0.
\end{equation*}
This means that when one of the equalities either $x=y$ or $y=z$ is satisfied, the function $g(x,y,z)$ possibly has a minimum. Let $x=y$, then the third equation from (\ref{ENC}) is fulfilled identically, and the first and second equations are identical and can be written as
\begin{equation*}
\displaystyle \ln \eta=\eta-1,~~\eta=\frac{z}{y}.
\end{equation*}
The last equation has only one root $\eta=1$, i.e. $y=z$.

Let us calculate the second derivatives and write the Hesse matrix:
\begin{equation}
G(x,y,z)=\left(\begin{array}{lll}
\displaystyle \frac{z-y}{x^2} & \displaystyle\frac{1}{x}-\frac{1}{y} & \displaystyle\frac{1}{z}-\frac{1}{x}\\
\displaystyle \frac{1}{x}-\frac{1}{y} & \displaystyle \frac{x-z}{y^2} & \displaystyle\frac{1}{y}-\frac{1}{z}\\
\displaystyle \frac{1}{z}-\frac{1}{x} & \displaystyle\frac{1}{y}-\frac{1}{z} & \displaystyle \frac{y-x}{z^2}
\end{array}
\right).
\end{equation}
Minors of the Hesse matrix are equal to
\begin{equation}
\begin{array}{c}
\displaystyle M_1(x,y,z)=\frac{z-y}{x^2},~M_2(x,y,z)=-\frac{(x-y)^2+(z-x)^2+(y-z)^2}{2x^2y^2},\\
\displaystyle M_3(x,y,z)=0.
\end{array}
\end{equation}

The Hesse matrix is not sign-defined, moreover, its determinant equals zero. Therefore, let us consider a small vicinity of the extremum point.

In a small vicinity $x=y=z$, provided that $\delta x\leq \delta y\leq \delta z$, the variation $\delta g$ of the function $g(x,y,z)$ is written in the form of
\begin{equation*}
\begin{array}{lll}
\delta g=(x+\delta y)\ln (x+\delta x)-(x+\delta x)\ln (x+\delta y)+(x+\delta x)\ln (x+\delta z)-\\
-(x+\delta z)\ln (x+\delta x)+(x+\delta z)\ln (x+\delta y)-(x+\delta y)\ln (x+\delta z)=\\
=(\delta z-\delta x)(\delta y-\delta x)(\delta z-\delta y)+o(((\delta x)^2+ (\delta y)^2+(\delta z)^2)^{3/2})\geq 0,
\end{array}
\end{equation*}
where values in the cubes of variations of the independent variables are nonzero, and the variation of the function $g(x,y,z)$ itself is positive by virtue of the Lemma conditions.  In the case when, for example, $\delta y=\delta x$, we have $g(x,y,z)\equiv 0$, and $f(x,y,z)\equiv 1$. Therefore, the extremum of the function is its non-strict minimum.
\end{proof}

\begin{figure}[H]
    \centering
    \includegraphics[width = 10cm]{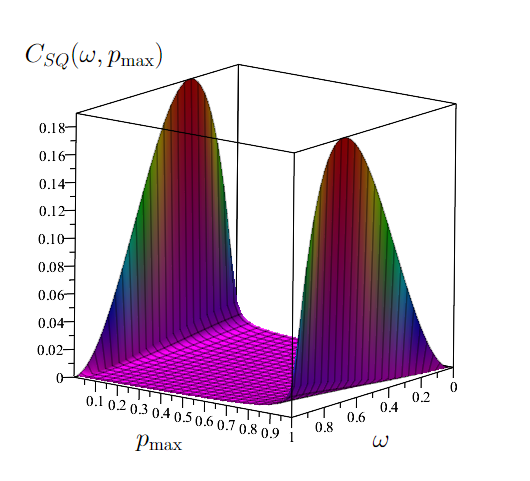}
    \caption{Level surfaces of statistical complexity $C_{SQ}(\omega,p_{\max})$}\label{Pic101}
\end{figure}

\begin{lemma}\label{lemmSQ}
The maximum statistical complexity (\ref{comp}) is achieved on the distribution of the form 
\begin{equation}\label{maximal}
    \left\{
\begin{array}{lll}
         \displaystyle p_i = \frac{1-p_{\max}}{N-1}, \quad i = \overline{1,N}~\backslash ~k,\\
         \displaystyle p_k = p_{\max},         
    \end{array}
    \right.
\end{equation}
where $p_{\max} = const$, i.e., at the appearance of a single component of an arbitrary index $k$ over the uniform distribution.
\end{lemma}
\begin{proof}
Without loss of generality, let us assume $k=N$. From equation (\ref{norm}) one variable $p_N$ from the set ${p_i}$ can be expressed through all the others:
\begin{equation}\label{p_n_equ}
    p_N = 1 - \sum_{i=1}^{N-1} p_i.
\end{equation}
Let us rewrite the equation \eqref{comp} in the form
\begin{equation}\label{compnew}
   C_{SQ}= -\frac{1}{\log N}\left(\sum_{i=1}^{N-1} p_i\log p_i + p_N \log p_N\right)\cdot \left(\sum_{i=1}^{N-1} \left(p_i - \frac{1}{N} \right)^2 + \left(p_N - \frac{1}{N}\right)^2\right).
\end{equation}
A necessary condition for the extremum of a function at an interior point of the domain (simplex \ref{discdefinit}) is that all partial derivatives of $p_i$ are equal to zero:
\begin{equation}\label{conditi}
    \frac{\partial C_{SQ}}{\partial p_i} = 0, \quad i  = 1, \dots, N - 1.
\end{equation}
Substituting the function \eqref{compnew} into \eqref{conditi} gives (provided that $\displaystyle \frac{\partial p_N}{\partial p_i} = -1$):
\begin{equation}
\begin{array}{ccc}
    \displaystyle \frac{\partial C_{SQ}}{\partial p_i} = -\frac{1}{\log N}\left( \log p_i -\log p_N  \right)\cdot \left(\sum_{i=1}^{N-1} \left(p_i - \frac{1}{N} \right)^2 + \left(p_N - \frac{1}{N}\right)^2\right) - \\
    \displaystyle - \frac{2}{\log N}\left(\sum_{i=1}^{N-1} p_i\log p_i + p_N \log p_N\right)\cdot  \left(p_i - p_{N} \right) = 0, \quad i  = 1, \dots, N - 1.
\end{array}    
\end{equation}
In a more convenient form the equations can be rewritten as
\begin{equation}\label{eq_proof_1}
    \displaystyle \frac{\partial C_{SQ}}{\partial p_i} = \frac{1}{\log N}\left( -\log p_i + \log p_N \right)\cdot D + 2H \cdot \left(p_i - p_N \right)= 0, \quad i  = 1, \dots, N - 1.
\end{equation}
Let us write the difference of any two equations from the system above for indices $i$ and $j$:
\begin{equation}\label{eq_proof_2}
\begin{array}{ccc}
    \displaystyle \frac{\partial C_{SQ}}{\partial p_i} - \frac{\partial C_{SQ}}{\partial p_j} = \frac{1}{\log N}\left( -\log p_i + \log p_j \right)\cdot D + 2H \cdot \left(p_i - p_j \right)= 0.
\end{array} 
\end{equation}
Given that the values of $D$ and $H$ are positive, the following equations can be constructed from the equations \eqref{eq_proof_1} and \eqref{eq_proof_2}, provided that the considered probabilities $p_j,~j=1,\ldots,N-1$ are not equal to $p_N$:
\begin{equation}\label{eq_proof_3}
\begin{array}{ccc}
    \displaystyle \frac{ \log p_i - \log p_j }{ \log p_N - \log p_j }-\frac{p_i - p_j }{p_N - p_j }=0,
\end{array} 
\end{equation}
\begin{equation}\label{eq_proof_4}
    \displaystyle (p_N-p_j)\log p_i+(p_i-p_N)\log p_j+(p_j-p_i)\log p_N=0,
\end{equation}
\begin{equation}\label{eq_proof_5}
    \displaystyle  p_i^{p_N-p_j} \cdot p_j^{p_i-p_N} \cdot p_N^{p_j-p_i}=1.
\end{equation}

After applying the Lemma \ref{lemma3} we conclude that the last equation can be satisfied when $p_i = p_j$.

Thus, it is obtained that each of the probabilities $p_i$ can take one of two different values, which define a distribution of the form
\begin{equation}\label{distrK}
\left\{
\begin{array}{l}
    \displaystyle p_i = \frac{1-p_{\max}}{K}, \quad \forall ~  i = 1, \dots, K,\\ 
 \displaystyle p_i = p_N = \frac{p_{\max}}{N-K}, \quad \forall ~  i = K+1, \dots, N.
\end{array}
\right.
\end{equation} 
Now we need to show that the maximum complexity corresponds to values $K=1$ and $K=N-1$.  For this purpose, let us calculate the value of the disequilibrium \eqref{disequ_1} on the distribution (\ref{distrK}), which we denote by $D^{(K)} (\omega,p_{\max})$:
\begin{equation}\label{disbKSQ}
D^{(K)}(\omega,p_{\max})=   \frac{1}{N}\frac{(p_{\max}+\omega-1)^2}{\omega(1-\omega)}, ~~~\omega=\frac{K}{N}.
\end{equation}
In turn, entropy is equal to
\begin{equation}\label{entrKSQ}
 \displaystyle H^{(K)}(\omega,p_{\max})= 1-\frac{1}{\log N}\left((1-p_{\max})\log\frac{1-p_{\max}}{\omega} +p_{\max}\log\frac{p_{\max}}{1-\omega}\right).
\end{equation}

The maximum of $C_{SQ}(\omega,p_{\max})$ at $N\leq 100$ was investigated numerically, and it was reached at $K=1$. From the expression for $D^{(K)}(\omega,p_{\max})$ (\ref{disbKSQ}), it can be seen that at $N\geq 101$ and when changing from $K=1$ to $K=2$ or from $K=N-1$ to $K=N-2$, its value changes by almost a factor of two, while the entropy (\ref{entrKSQ}) changes only slightly.
Thus, the probability distribution (\ref{distrK}) that delivers the complexity function to the extremum value at $K=1$ or $K=N-1$ is of the form (\ref{maximal}).
\end{proof}
For clarity Fig. \ref{Pic101} shows the graph $C_{SQ}=C_{SQ}(\omega,p_{\max})$ at $N=1024$, where $\omega$ is changing continuously (although $K$ is changing discretely).

\begin{corollary}
Let us substitute the values of $p_i$ and $p_N=p_{\max}$ from (\ref{maximal}) into (\ref{compnew}) and consider the complexity $C_{SQ}$ as a function of $p_{\max}$. For sufficiently large values of $N$, it will take the following form
$$
C_{SQ} \approx (1-p_{\max})\cdot p_{\max}^2.
$$    
Whence it follows that this function takes the maximum value $C_{SQ}^* \approx 4/27$ when $p_{\max}=2/3$.
\end{corollary}

\begin{corollary}\label{corSQmin}
The minimum value of $C_{SQ}=0$ is achieved on a uniform distribution $\displaystyle p = (1/N, \dots, 1/N)$.
\end{corollary}

\begin{figure}[H]
    \centering
    \includegraphics[width = 10cm]{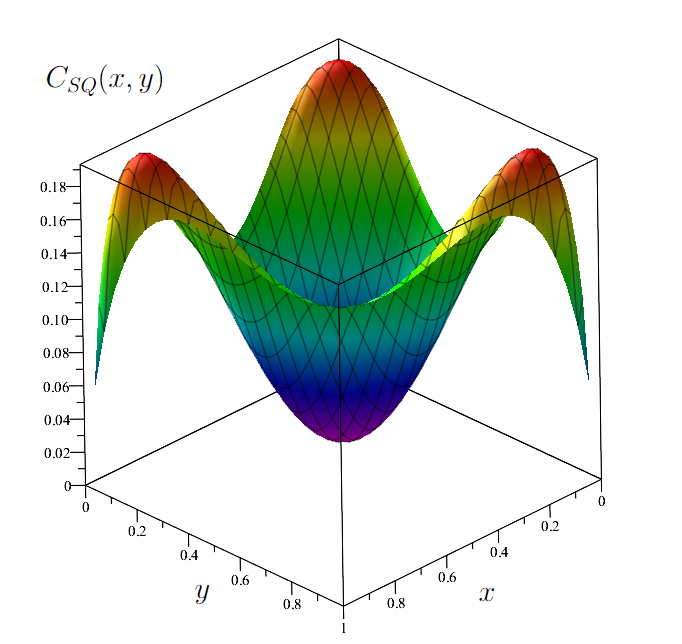}
    \caption{Level surfaces of statistical complexity $C_{SQ}(x,y)$ for $p=\{p_1=x,~p_2=y,~p_3=1-x-y\}$}\label{Pic10}
\end{figure}

The validity of the Lemma \ref{lemmSQ} for the case when the discrete distribution $p=\{p_1,p_2,p_3\}$ consists of three samples is demonstrated in Fig. \ref{Pic10}. The complexity depends on two variables, since one of the probabilities can be expressed through the others. Here $C_{SQ}$ has three identical pronounced maxima and three identical local minima pertaining to the cases $p_1=p_2$, $p_2=p_3$, $p_1=p_3$ when the necessary extremum conditions are met, and a global minimum when $p_1=p_2=p_3$.

Table \ref{tableSQ} shows the change of optimal parameters $C_{SQ}(w,p_{\max})$ with increasing $N$. 
\def\arraystretch{1.8}
\begin{table}[h]
\caption{Optimal parameters $C_{SQ}(\omega, p_{max})$ for different values of $N$.}\label{tableSQ}
\centering
\begin{tabularx}{0.88\textwidth} { 
  | >{\centering\arraybackslash}X 
  | >{\centering\arraybackslash}X 
  | >{\centering\arraybackslash}X 
  | >{\centering\arraybackslash}X 
  | >{\centering\arraybackslash}X | }
 \hline
{ $N$ } & $C_{SQ}(\omega^*,p_{\max}^*)$ & $p_{\max}^*$  & $\omega^*$ & $N-K^*$ \\ 

 \hline
$3$  & $0,1932$ & $0,8315$  & $0,6666$  & $1$ \\
 \hline
$256$  & $0,1994$ & $0,7044$  & $0,9960$  & $1$ \\

 \hline
$512$  & $0,1942$ & $0,7008$  & $0,9980$  & $1$ \\
 \hline
 $1024$  & $0,1898$ & $0,6979$  & $0,9990$  & $1$ \\
 \hline
  $2048$  & $0,1861$ & $0,6955$  & $0,9995$  & $1$ \\
 \hline
\end{tabularx}
\end{table}

The necessary extremum conditions $C_{SQ}$ for the discrete distribution $p=\{p_1=x,~p_2=y,~p_3=1-x-y\}$ are written out according to \eqref{eq_proof_1} as follows:
\begin{equation}\label{necess_cond}
\begin{cases}
    \displaystyle \left( -\log x + \log (1-x-y) \right) \left(\left(x-\frac{1}{3}\right)^2+\left(y-\frac{1}{3}\right)^2+\left(1-x-y-\frac{1}{3}\right)^2\right) -\\- 2(x\log x+y\log y+(1-x-y)\log (1-x-y))  \left( -1+y) \right)= 0,\\
    \displaystyle\left( -\log y + \log (1-x-y) \right) \left(\left(x-\frac{1}{3}\right)^2+\left(y-\frac{1}{3}\right)^2+\left(1-x-y-\frac{1}{3}\right)^2\right) -\\- 2(x\log x+y\log y+(1-x-y)\log (1-x-y))  \left(-1+x \right)= 0.\\
    \end{cases}
\end{equation}

\begin{figure}[H]
    \centering
    \includegraphics[width = 14cm]{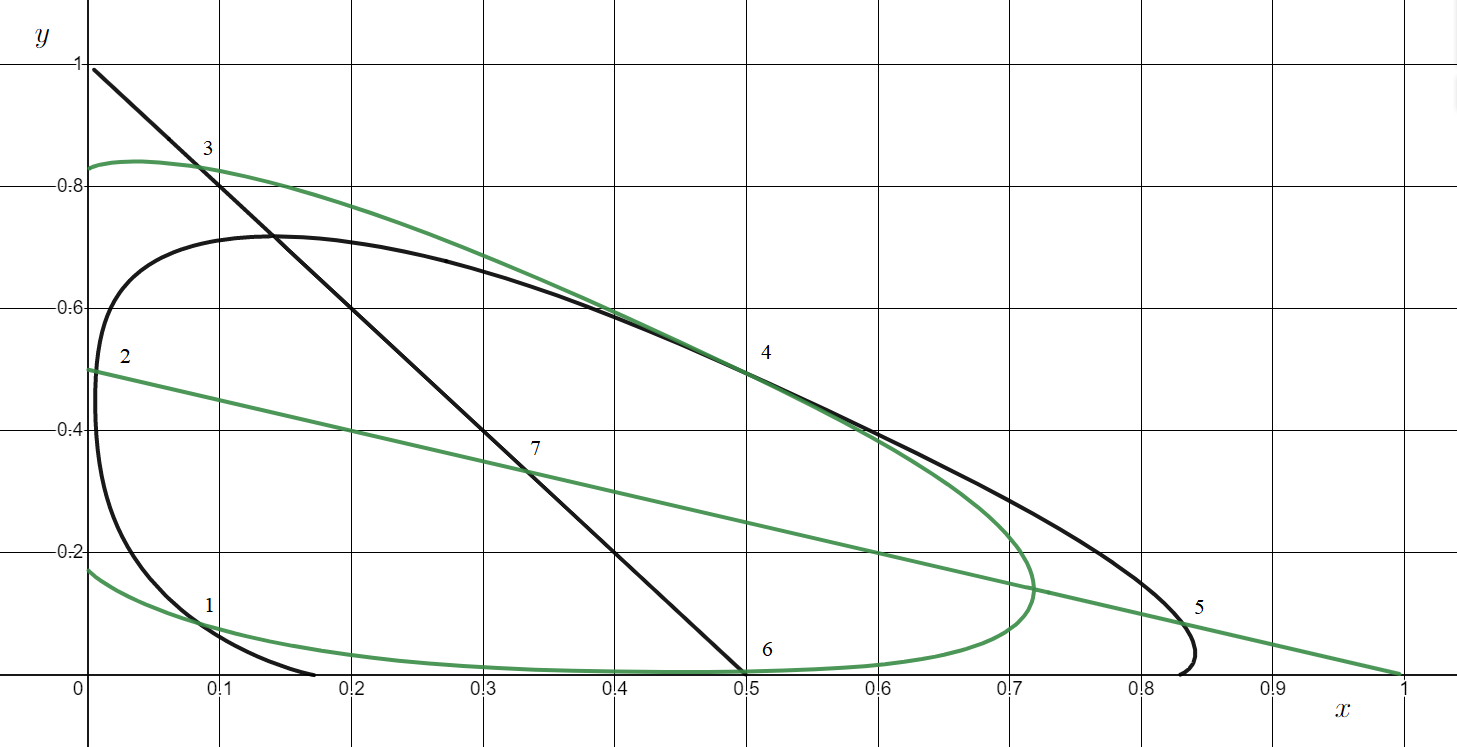}
    \caption{Curves corresponding to the necessary conditions of extremum $C_{SQ}$ for $p=\{p_1=x,p_2=y,p_3=1-x-y\}$}\label{extr_cond_3}
\end{figure}

The implicit equations of the system \eqref{necess_cond} describe the curves shown in Fig. \ref{extr_cond_3}.

The black and green curves correspond to the first and second implicit equations of the \eqref{necess_cond} system, respectively. In Fig. \ref{extr_cond_3} seven extremum points are marked, for which the value of statistical complexity is calculated. All the obtained data is summarized in Table \ref{tableextrem}.

\def\arraystretch{1.8}
\begin{table}[h]
\caption{Extremum points of statistical complexity at $N=3$.}\label{tableextrem}
\centering
\begin{tabularx}{0.9\textwidth} { 
  | >{\centering\arraybackslash}X 
  | >{\centering\arraybackslash}X 
  | >{\centering\arraybackslash}X 
  | >{\centering\arraybackslash}X 
  | >{\centering\arraybackslash}X 
  | >{\centering\arraybackslash}X 
  | >{\centering\arraybackslash}X 
  | >{\centering\arraybackslash}X | }
 \hline
{ $p$ } & 1 & 2  & 3 & 4 & 5 & 6 & 7\\ 

 \hline
$p_1$&   $0,08425$ & $0,006$  & $0,08425$  & $0,497$ & $0,8315$ & $0,497$ & $0,(3)$\\
 \hline
$p_2$ & $0,08425$ & $0,497$  & $0,8315$  & $0,497$ & $0,08425$ & $0,006$ & $0,(3)$\\
 \hline
$p_3$ & $0,8315$ & $0,497$  & $0,08425$  & $0,006$ & $0,08425$ & $0,497$&$0.(3)$\\
\hline
$C_{SQ}$ & $0,1932$ & $0,1062$  & $0,1932$  & $0,1062$ & $0,1932$ & $0,1062$&$0$\\
 \hline
\end{tabularx}
\end{table}

The first, third and fifth points of maximum correspond to the same value of the maximum of the function. The second, fourth and sixth minimum points correspond to the same value of the minimum of the function. It is worth noting that the extremum points correspond to the case $K=N-1=2$ except for the global minimum, where all probabilities are equal to each other, and thus describe three local minima, one global minimum, and three equal maxima of statistical complexity in Fig. \ref{Pic10}. We can separately plot the statistical complexity in the case $p=\{p_1=x,p_2=x,p_3=1-2x\}$.

\begin{figure}[H]
    \centering
    \includegraphics[width = 14cm]{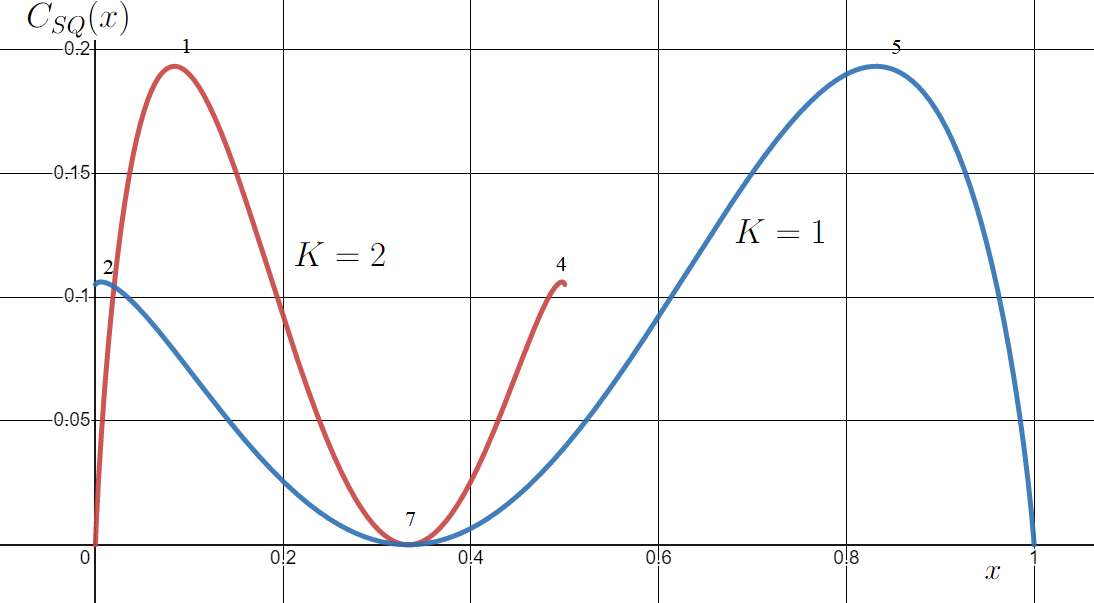}
    \caption{Statistical complexity $C_{SQ}$ for $p=\{p_1=x,~p_2=x,~p_3=1-2x\}$}\label{csq_2x}
\end{figure}
In Fig. \ref{csq_2x} the extremum points are marked according to Table {\ref{tableextrem}, which cover all cases $p=\{p_1=x,~p_2=y,~p_3=1-x-y\}$.

\subsection{Optimization of $C_{JSD}$}\label{JSDoptimize}

Let us apply a similar approach using the Jensen-Shannon divergence as the disequilibrium to the statistical complexity $C_{JSD}$
\begin{equation}
    C_{JSD}(p) = H(p)\cdot JSD(p||q), \quad q_j=1/N,~j=1,\ldots,N,
\end{equation}
which can be written by considering the expression $JSD(p||q)$ through entropy
\begin{equation}\label{CJSD}
    C_{JSD}(p) = H(p)\cdot \left(H(m) - \frac{1}{2}(H(p)+H(q))\right)\cdot\log N, \quad m=\frac{p+q}{2}.
\end{equation}

It is possible to write out the necessary conditions of extremum for the statistical complexity of the form (\ref{CJSD}), but at the same time a Lemma similar to the Lemma \ref{lemmSQ} cannot be proved.

The equation (\ref{CJSD}) is written in variables $p_i$ and $p_N$ in accordance with the approach of the Lemma \ref{lemmSQ}
\begin{equation}
\hspace{-1cm}
    C_{JSD}(p) = -\left(\sum_{i=1}^{N-1} p_i\log p_i + p_N \log p_N\right)\cdot \left(H(m) - \frac{1}{2}\left(-\frac{1}{\log N}\left(\sum_{i=1}^{N-1} p_i\log p_i + p_N \log p_N\right)+1\right)\right),
\end{equation}
where
\begin{equation}
    \displaystyle H(m) = -\frac{1}{\log N}\left(\sum_{i=1}^{N-1} \frac{p_i + \frac{1}{N}}{2}\log \frac{p_i + \frac{1}{N}}{2} + \frac{p_N + \frac{1}{N}}{2}\log \frac{p_N + \frac{1}{N}}{2}\right).
\end{equation}
Then taking into account \eqref{p_n_equ}
\begin{equation}
    \displaystyle \frac{\partial H(m)}{\partial p_i} = -\frac{1}{\log N}\left(\frac{1}{2}\log \frac{p_i + \frac{1}{N}}{2} - \frac{1}{2}\log \frac{p_N + \frac{1}{N}}{2}\right),~i=1,\ldots,N-1.
\end{equation}
The necessary conditions of extremum are obtained in the following form after combining all partial derivatives:
\begin{equation}
\begin{array}{cc}
    \displaystyle \frac{\partial C_{JSD}(p)}{\partial p_i} =\displaystyle \frac{H(p)}{2}\cdot \left(-\left(\log \frac{p_i + \frac{1}{N}}{2} - \log \frac{p_N + \frac{1}{N}}{2}\right) +\left( \log p_i -\log p_N  \right)\right)- \\
    \displaystyle-\left( \log p_i -\log p_N  \right)\cdot JSD(p||q)= 0,~i=1,\ldots,N-1.
\end{array}
\end{equation}
The equations after simplification are given:
\begin{equation}\label{NECCJSD}
\begin{array}{cc}
    \displaystyle \displaystyle H(p)\cdot \left(\log \frac{p_i + \frac{1}{N}}{2} - \log \frac{p_N + \frac{1}{N}}{2}  \right)
    +\left( \log p_i -\log p_N  \right)\cdot \left(2JSD(p||q)-H(p)\right) = 0,~i=1,\ldots,N-1.
\end{array}
\end{equation}
Then the difference of equations (\ref{NECCJSD}) for indices $i,j$ takes the following form
\begin{equation}\label{eq75}
\begin{array}{cc}
    \displaystyle \left( \log p_i -\log p_j  \right)\cdot \left(2JSD(p||q)-H(p)\right) + H(p)\cdot \left(\log \frac{p_i + \frac{1}{N}}{2} - \log \frac{p_j + \frac{1}{N}}{2}  \right) = 0.
\end{array}
\end{equation}

\begin{remark}\label{remjsd}
    It follows from the form of the system of equations (\ref{eq75}) that the system is satisfied if $p_i=p_j$, which is one of the necessary conditions for the extremum of the function (\ref{CJSD}).
    Due to the nonlinearity of the system consisting of equations (\ref{eq75}), it may have other roots.
\end{remark}

Fig. \ref{Pic11} shows a surface plot of the statistical complexity level of the form (\ref{CJSD}) when the discrete distribution $p=\{p_1,~p_2,~p_3\}$ consists of three samples to illustrate the Remark \ref{remjsd}.

\begin{figure}[H]
    \centering
    \hspace{-1cm}
    \includegraphics[width = 9.5cm]{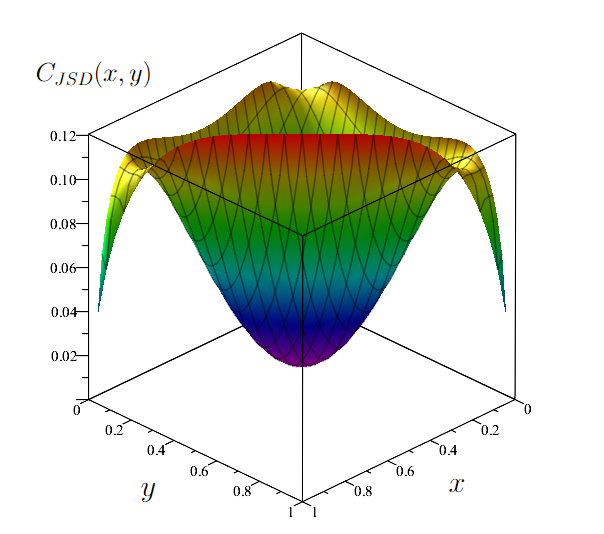}
        \caption{Level surfaces of statistical complexity $C_{JSD}$ for $p=\{p_1=x,p_2=y,~p_3=1-x-y\}$}\label{Pic11}
\end{figure}

It can be seen from Fig. \ref{Pic11} that the points satisfying $p_1=p_2$, $p_2=p_3$ and $p_2=p_3$ are the saddle points of the surface if the necessary extremum conditions are satisfied.

It was previously established that the distribution \eqref{distrK} delivers an extremum to $C_{SQ}$ at $K=N-1$. Next, it will be shown that it also delivers the extremum to the complexity based on the total variation of the measure $TV(p,q)$. Therefore, we propose to find the maximum of $C_{JSD}$ on this distribution and compare the obtained optimal distribution parameters at fixed $N$. Let us write out the complexity explicitly and obtain
\begin{equation}\label{CJSD_k_p}
    C^{(K)}_{JSD} = H^{(K)}\cdot \left(H^{(K)}(m) - \frac{1}{2}(H^{(K)}+1)\right)\cdot \log N,
\end{equation}
where $H^{(K)}$ is corresponding to \eqref{entrKSQ}, and $H^{(K)}(m)$ is given by the following formula:
\begin{equation}
 \displaystyle H^{(K)}(m)= 1-\frac{1}{\log N}\left(\frac{(1-p_{\max}+\omega)}{2}\log\frac{(1-p_{\max}+\omega)}{2\omega} +\frac{(1+p_{\max}-\omega)}{2}\log\frac{(1+p_{\max}-\omega)}{2-2\omega}\right).
\end{equation}

Table \ref{tableJSD} shows the change of optimal parameters $C_{JSD}$ with the growth of $N$. 

\def\arraystretch{1.8}
\begin{table}[h]
\caption{Optimal parameters $C_{JSD}(\omega,p_{\max})$ for different values of $N$.}\label{tableJSD}
\centering
\begin{tabularx}{0.99\textwidth} { 
  | >{\centering\arraybackslash}X 
  | >{\centering\arraybackslash}X 
  | >{\centering\arraybackslash}X 
  | >{\centering\arraybackslash}X 
  | >{\centering\arraybackslash}X | }
 \hline
{ $N$ } & $C_{JSD}(\omega^*,p_{\max}^*)$ & $p_{\max}^*$  & $\omega^*$ & $N-K^*$ \\ 
 \hline
$3$  & $0,1266$ & $1$  & $0,4083$  & $1\text{ or }2$ \\
 \hline
$256$  & $0,4482$ & $1$  & $0,8703$  & $33$ \\
 \hline
$512$  & $0,4790$ & $1$  & $0,8897$  & $56$ \\
 \hline
 $1024$  & $0,5065$ & $1$  & $0,9051$  & $97$ \\
 \hline
  $2048$  & $0,5312$ & $1$  & $0,9171$  & $170$ \\
 \hline
\end{tabularx}
\end{table}

For clarity, Fig. \ref{PicJSD_FIG} shows the graph
$C_{JSD}=C_{JSD}(\omega,p_{\max})$ at $N=1024$, where $\omega$ is changing continuously (although $K$ is changing discretely). 

\begin{figure}[H]
    \centering
    \hspace{-1cm}
    \includegraphics[width = 8cm]{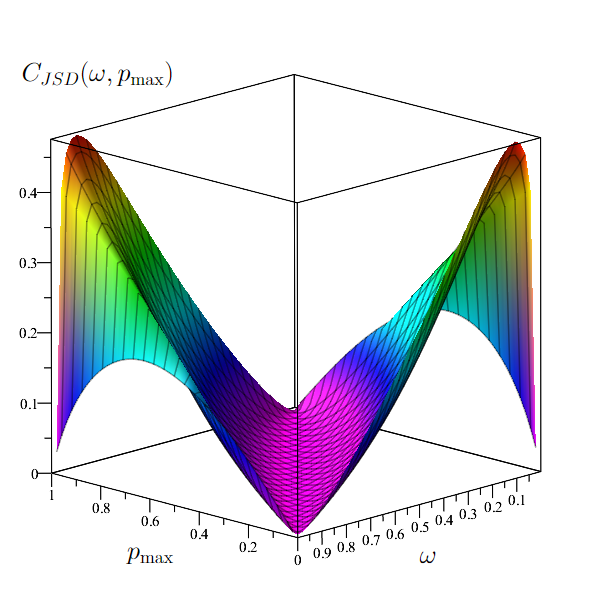}
    \caption{Level surfaces of statistical complexity $C_{JSD}(\omega,p_{\max})$}\label{PicJSD_FIG}
\end{figure}

The results shown in Table \ref{tableJSD} demonstrate that for the chosen class of distributions \eqref{distrK} the set of $N$ components, where $K$ is equal to each other and the rest are zero, is optimal. It is worth noting that for the resulting distribution $C_{JSD}$ is not zero, as well due to the summand $H^{(K)}(m)$, which corresponds to the already "shifted" $~$ distribution consisting of $K$ elements equal to $\displaystyle\frac{\frac{1}{K}+\frac{1}{N}}{2}$, and $N-K$ samples of $\displaystyle\frac{1}{2N}$ each.

\subsection{Optimization $C_{TV}$}
Let us proceed to analyze statistical complexity based on total variation 
  \begin{equation}\label{CTV}
        \displaystyle C_{TV}(p)=-\frac{1}{4\log N}\left(\sum_{i=1}^N p_i\log p_i \right)\cdot\left(\sum _{i=1}^N\left|p_i-\frac{1}{N}\right|\right)^2.
    \end{equation}
\begin{proposition}
According to the expression for the error function \eqref{ErrTVR} from the Neyman-Pearson Lemma \ref{NeymanPearson} and the definition of \eqref{CTV}, we propose to use $C_{TV}$ as a criterion to solve the problem of hypothesis testing and indicating the appearance of a deterministic component of a useful signal in noise.
\end{proposition}
The following Lemma is valid.
 \begin{lemma}
The maximum statistical complexity (\ref{CTV}) is achieved on the family of distributions (\ref{distrK}).
 \end{lemma}   
\begin{proof}    

Given the symmetry of the function \eqref{CTV} and the simplex \eqref{discdefinit}, without restriction of generality, we find an integer $K \in\{1,...,N-1\}$ for which the maximum of this Lemma is achieved on the part of the simplex \eqref{discdefinit} defined by the constraints $p_i\leq 1/N$ for $i=1,...,K$ and $p_i\geq 1/N$ for $i=K+1,...,N$. Let us rewrite the equation (\ref{CTV}) in the form of
\begin{equation}
   C_{TV}= -\frac{1}{4\log N}\left(\sum_{i=1}^{N-1} p_i\log p_i + p_N \log p_N\right)\cdot \left(\sum_{i=1}^{K} \left(-p_i + \frac{1}{N} \right) + \sum_{i=K+1}^{N}\left(p_i - \frac{1}{N}\right)\right)^2.
\end{equation}
Then for $i=1,\ldots,K$ the necessary conditions of extremum take the form
\begin{equation}\label{NECCTV1}
\begin{array}{ccc}
    \displaystyle \frac{\partial C_{TV}}{\partial p_i} = -\frac{1}{\log N}\left( \log p_i -\log p_N  \right)\cdot D_{TV} - 2H(p)\sqrt{D_{TV}}=0,
     \quad i  = 1, \dots, K,
\end{array}    
\end{equation}
and for $i=K+1,\ldots,N$ the following is true 
\begin{equation}\label{NECCTV2}
\begin{array}{ccc}
    \displaystyle \frac{\partial C_{TV}}{\partial p_i} = -\frac{1}{\log N}\left( \log p_i -\log p_N  \right)\cdot D_{TV} =0,
     \quad i  = K+1, \dots, N.
\end{array}    
\end{equation}
Let us compose the difference of two equations from (\ref{NECCTV1}) for indices $i$ and $j$. Whence it follows that if $D_{TV}\not=0$, then $p_i=p_j$ at $i=1,\ldots,K$. Whereas it follows from (\ref{NECCTV2}) that $p_i=p_N$ at $i=K+1,\ldots,N$. Again we obtain that the family of distributions (\ref{distrK}) delivers the maximum of the complexity function, now $C_{TV}$.
\end{proof}

Next, we need to determine the optimal values of $K$ and $p_{\max}$. 
For this purpose, let us calculate the disequilibrium value $D^{(K)}_{TV}$ on the distribution (\ref{distrK}): 
\begin{equation}\label{DKTV}
D^{(K)}_{TV}=  (p_{\max}+\omega-1)^2, ~~~\omega=\frac{K}{N}.
\end{equation}
In turn, entropy is equal to
\begin{equation}\label{HKTV}
 \displaystyle H^{(K)}= 1-\frac{1}{\log N}\left((1-p_{\max})\log\frac{1-p_{\max}}{\omega} +p_{\max}\log\frac{p_{\max}}{1-\omega}\right).
\end{equation}

For clarity, Fig. \ref{Pic13} shows the graph $C_{TV}=C_{TV}(\omega,~p_{\max})$ at $N=1024$, where $\omega$ is changing continuously (although $K$ is changing discretely).

\begin{figure}[H]
    \centering
    \hspace{-1cm}
    \includegraphics[width = 9.5cm]{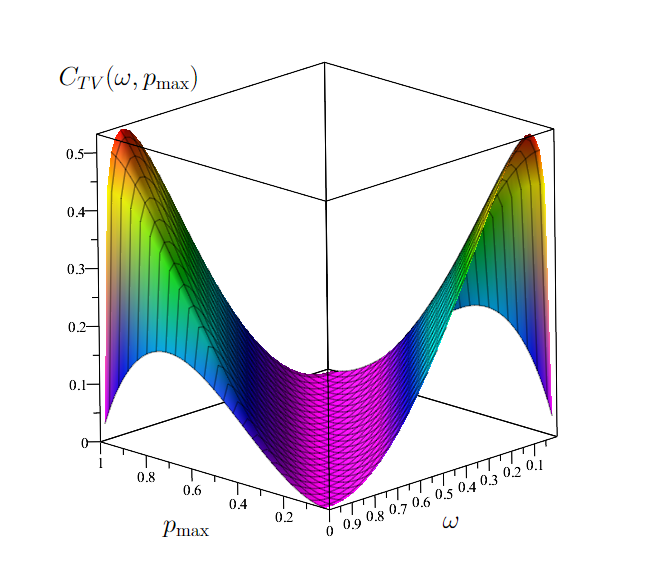}
    \caption{Level surfaces of statistical complexity $C_{TV}(\omega,p_{\max})$}\label{Pic13}
\end{figure}

Let us compose the necessary conditions \eqref{NECCTV1} and \eqref{NECCTV2} of the extremum of the statistical complexity $C_{TV}$ written out through the variables $p_{\max},~\omega$.
\small
\begin{equation}\label{nece_cond_ctv}
    \begin{cases}
    \displaystyle f_{1}^N(p_{\max},\omega) \coloneqq 2(p_{\max}+\omega-1)\left(\left( 1-\frac{(1-p_{\max})\log \frac{1-p_{\max}}{\omega}+p_{\max}\log\frac{p_{\max}}{1-\omega} }{\log N}  \right)\right.-\\
    \displaystyle\left.-\frac{(p_{\max}+\omega-1)}{2\log N}\left(\log\frac{p_{\max}}{1-\omega}-\log\frac{1-p_{\max}}{\omega}\right)\right)=0,\\\\
    
    \displaystyle f_{2}^N(p_{\max},\omega) \coloneqq 2(p_{\max}+\omega-1)\left(\left( 1-\frac{(1-p_{\max})\log\frac{1-p_{\max}}{\omega}+p_{\max}\log\frac{p_{\max}}{1-\omega}}{\log N}  \right)-\right.\\
\displaystyle\left.-\frac{(p_{\max}+\omega-1)}{2\log N}\left(\frac{p_{\max}}{1-\omega}-\frac{1-p_{\max}}{\omega}\right)\right)=0.
    \end{cases}
\end{equation}
\normalsize
The intersections of the curves corresponding to the implicit equations \eqref{nece_cond_ctv} are related to to the extremum points of $C_{TV}$. Let us compose the difference of two necessary conditions of the extremum
\begin{equation}\label{f1-f2}
\begin{array}{c}
    \displaystyle f_{3}^N(p_{\max},\omega) \coloneqq \displaystyle f_{1}^N(p_{\max},\omega)-\displaystyle f_{2}^N(p_{\max},\omega) = \\
    \displaystyle=\frac{(p_{\max}+\omega-1)^2}{\log N}\left(-\log\frac{p_{\max}}{1-\omega}+\log\frac{1-p_{\max}}{\omega}+\frac{p_{\max}}{1-\omega}-\frac{1-p_{\max}}{\omega}\right)=0.
    \end{array}
\end{equation}

Let us construct the implicit curves of equations $f_{1}^N(p_{\max},\omega),~f_{2}^N(p_{\max},\omega),~f_{3}^N(p_{\max},\omega)$ for some values of $N$. For convenience, the index $N$ of $f_{3}^N(p_{\max},\omega)$ can be omitted since the implicit curve of the equation \eqref{f1-f2} is independent of $N$.

\begin{figure}[H]
    \centering
    \hspace{-1cm}
    \includegraphics[width = 9.5cm]{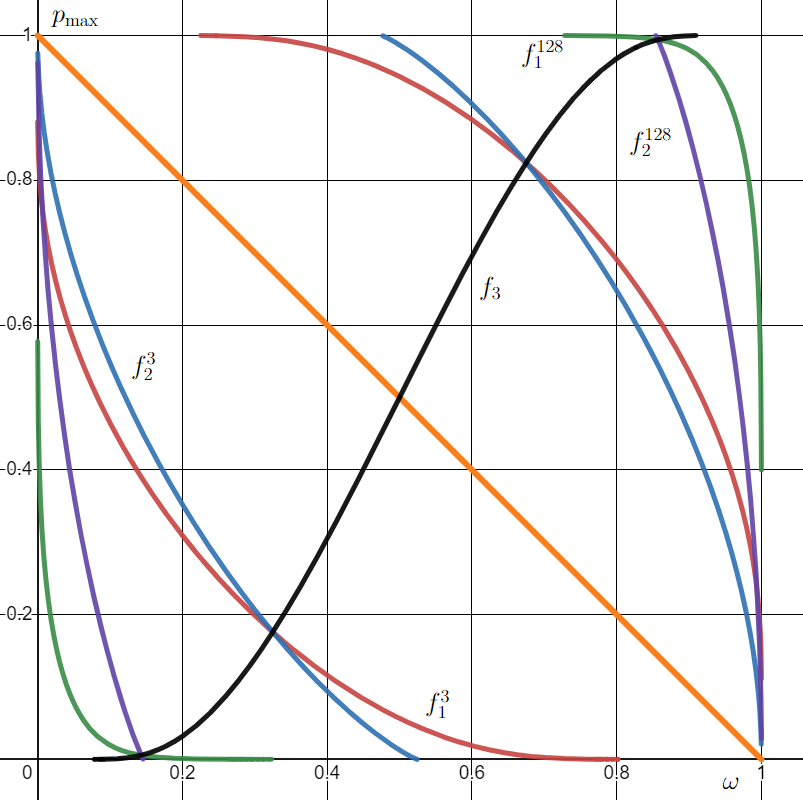}
    \caption{Curves of necessary conditions of extremum $C_{TV}$ for $N=3$ and $N=128$. Equations $f_{1}^N(p_{\max},\omega)$ and $f_{2}^N(p_{\max},\omega)$ have common solution $(p_{\max}+\omega-1=0)$ (orange curve)}\label{CTV_pic}
\end{figure}

According to Fig. \ref{CTV_pic} the statistical complexity has, in addition to the minimum points $(p_{\max}+\omega-1=0=0)$, where $C_{TV}=0$, also two maximum points for each value of $N$: $(p_{\max}^*,\omega^*)$ and $(1-p_{\max}^*,1-\omega^*)$, which lie on the curve $f_{3}(p_{\max},\omega)$.

In Table \ref{tableTV} the optimal values of the parameters of the formulas (\ref{DKTV}) and (\ref{HKTV}) that maximize the statistical complexity of $C_{TV}$ are given.

\def\arraystretch{1.8}
\begin{table}[h]
\caption{Optimal parameters $C_{TV}(\omega,p_{\max})$ for different values of $N$.}\label{tableTV}
\centering
\begin{tabularx}{0.9\textwidth} { 
  | >{\centering\arraybackslash}X 
  | >{\centering\arraybackslash}X 
  | >{\centering\arraybackslash}X 
  | >{\centering\arraybackslash}X 
  | >{\centering\arraybackslash}X | }
 \hline
{ $N$ } & $C_{TV}(\omega^*,p_{\max}^*)$ & $p_{\max}^*$  & $\omega^*$ & $N-K^*$ \\ 

 \hline
$3$  & $0,1289$ & $0,8241$  & $0,6751$  & $1$ or $2$\\
 \hline
$256$  & $0,4789$ & $0,9976$  & $0,8752$  & $32$ \\
 \hline
$512$  & $0,5120$ & $0,9991$  & $0,8901$  & $56$ \\
 \hline
 $1024$  & $0,5410$ & $0,9997$  & $0,9022$  & $100$ \\
 \hline
 $2048$  & $0,5667$ & $0,9999$  & $0,9122$  & $180$ \\
  \hline
\end{tabularx}
\end{table}

Additionally, the case $N=3$ is shown in Fig. \ref{Pic14}, which shows a graph of the level surface $C_{TV}$ when the discrete distribution $p=\{p_1,p_2,p_3\}$ consists of three samples.

\begin{figure}[H]
    \centering
    \hspace{-1cm}
    \includegraphics[width = 10.5cm]{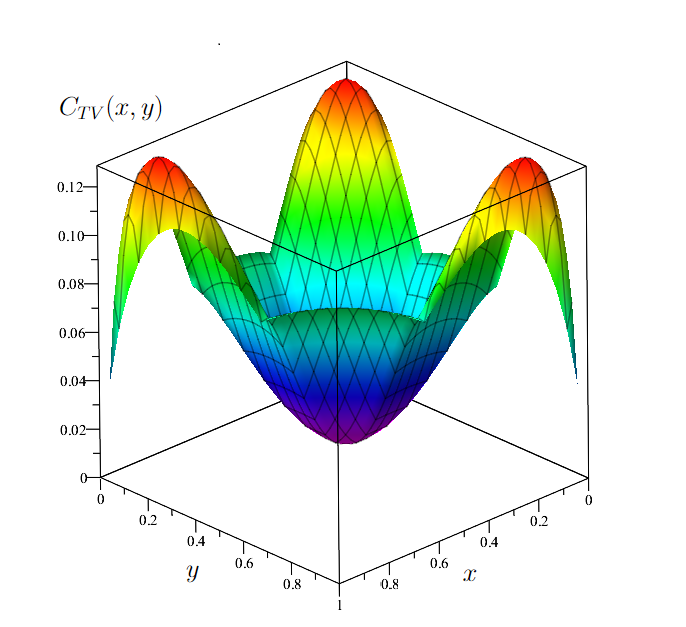}
    \caption{Level surfaces of statistical complexity $C_{TV}(x,y)$ for $p=\{p_1=x,~p_2=y,~p_3=1-x-y\}$}\label{Pic14}
\end{figure}

\section{Statistical Complexity Modeling and Comparison}\label{Modelling}
We analyze the optimal parameters that maximize different types of statistical complexity and compare the values in Tables. \ref{tableSQ}, \ref{tableJSD}, and \ref{tableTV}. Of main interest are the maximum complexity values and the optimal values of $K$. The maximum values of $C_{TV}(\omega^*,p_{\max}^*)\in[0,1]$, $C_{JSD}(\omega^*,p_{\max}^*)\in[0,1]$ are close to each other and grow with increasing $N$. The optimal values of $K$ for these two types of complexity are also close.

To demonstrate the analytical results obtained in the previous sections, the application of three variants of statistical complexity in the problem of useful signal indication in a noise mixture for synthesized signals is shown. An algorithm from \cite{Sensors} based on the computation of discrete distributions $p$ from the spectral representation of time series is applied.

The synthesized 10-second signal is the sum of a finite number of cosine oscillations mixed with white noise:
\begin{equation}\label{signal_mod1}
\displaystyle x(t)=I(t)\sum_{i=1}^K A_i\cos(2\pi f_i t + \Delta \phi_i) + w(t), \quad t \in [0, 10],\\
\end{equation}

where $A_i, ~f_i,~\Delta \phi_i$ are the amplitudes, frequencies, and random phases of the harmonic oscillations, respectively, $w(t)$ is the white noise, and $I(t)$ is the indicator function for the presence of the useful signal in the signal-noise mixture.

$I(t)$ is chosen so that the harmonic signals are present in the middle of the final sequence $x(t)$.

\begin{equation}
    I(t) = \left\{
    \begin{array}{lll}
        0, & t \in [0, 3), \\
        1, & t \in [3, 7], \\
        0, & t \in (7, 10].
    \end{array}
    \right.
\end{equation}

The algorithm has the following structure:
\begin{enumerate}
    \item The signal synthesized with sampling frequency $f_s$ is divided into short windows containing $N = 2048$ samples each.
    \item Next, the spectrum for each window is calculated using the FFT algorithm.
    \item Based on the spectrum, the discrete densities $p_i,~i = 1, \dots, N$ are calculated by normalizing it.
    \item The information characteristics $C_{SQ}(p),~C_{JSD}(p),~C_{TV}(p)$ are calculated for the obtained set $p_i$.
    \item The obtained sequence of information characteristic values is shown along with the signal on the time axis.
\end{enumerate}

It should be noted that the parameters $f_s$ and $N$ were chosen to exclude the effect of spectrum spreading, i.e., to obtain clear spectral components corresponding to $K$ harmonic functions from the formula \eqref{signal_mod1}. The signal-to-noise ration is chosen to be close to one.

The threshold $\gamma$ for the decision rule is proposed to be chosen as $25\%$ of the maximum criterion value for selected $N$ from the tables \ref{tableSQ}, \ref{tableJSD}, \ref{tableTV}:
\begin{equation}
    \begin{array}{ccc}
        \gamma_{CQ} = 0,25 \cdot 0,1861 = 0,0465; \\
        \gamma_{JSD} = 0,25 \cdot 0,5312 = 0,1328; \\
        \gamma_{TV} = 0,25 \cdot 0,5667 = 0,1417.
    \end{array}
\end{equation}

The convenience of choosing such a threshold is that it does not depend on a particular noise realization and is based on analytically derived maximum values of statistical complexity functions.

In all plots, the blue color indicates the amplitude of the original signal, and the red color indicates the statistical complexity, which is calculated using the algorithm described above. The horizontal axis represents time in seconds, and the vertical axes represent the magnitude of the signal amplitude (left) and the criterion (right). The black dashed line shows the value of the selected threshold $\gamma$.

In the first experiment, the number of sinusoidal signals and respectively spectral components is equal to $K = 3$ at $N=2048$. Fig. \ref{firstexp} shows the dependencies of statistical complexities on time for the synthesized signal.

As can be seen, the values of $C_{SQ}$ and $C_{TV}$ exceed the selected threshold for the interval of signal presence, which allows us to confidently conclude that the signal has occurred. As for $C_{JSD}$, the a priori threshold selection was unsuccessful because the true value of its maximum is unknown, as shown in subsection \ref{JSDoptimize}. If we change the threshold value upwards by $20\%$, the detection based on $C_{JSD}$ will be as successful as that based on $C_{TV}$.

\begin{figure}[h]
    \centering
    \includegraphics[width = 12cm]{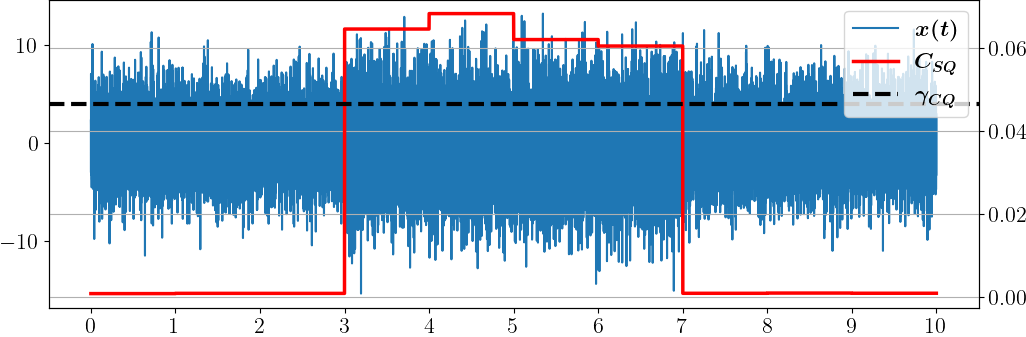}
    \vspace{0.1cm}
    \includegraphics[width = 12cm]{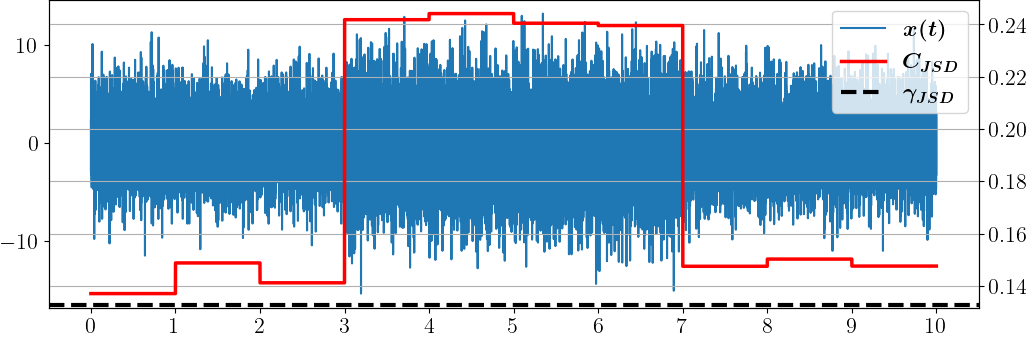}
    \vspace{0.1cm}
    \includegraphics[width = 12cm]{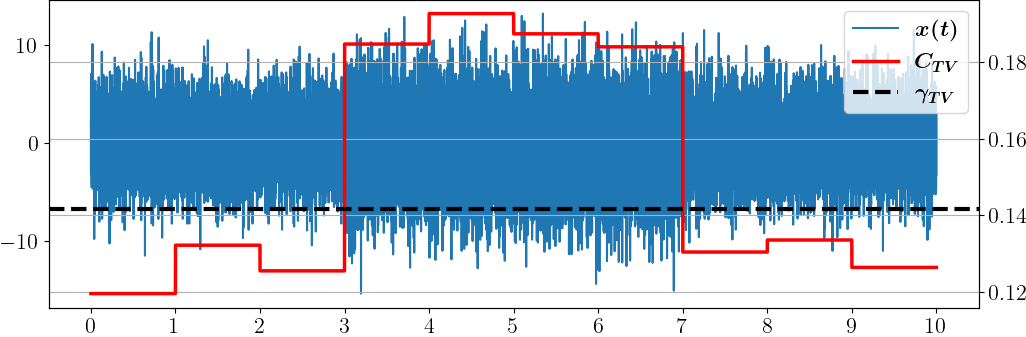}
    \caption{Three components, $K = 3$}
    \label{firstexp}
\end{figure}

In the second experiment, the number of spectral components $K = 30$. In this case, $C_{SQ}$ ceases to show a satisfactory result in the sense of exceeding the chosen threshold, since the function $C_{SQ}$ degrades strongly with increasing of $K$, but still allows a signal indication, as can be seen in Fig. \ref{secondexp}. The complexity function $C_{TV}$ still confidently exceeds the threshold, as in the first experiment, and $C_{JSD}$ exceeds the threshold on the whole signal, as in the first experiment. 

Thus, we can conclude that $C_{TV}$ is the most convenient in a practical sense, since it works well on signals with a large number of spectral components and allows one to decide on the appearance of an useful signal, using a fairly simple rule, associated with the choice of threshold based on the theoretical maximum value for the statistical complexity.

\begin{figure}[H]
    \centering
    \includegraphics[width = 12cm]{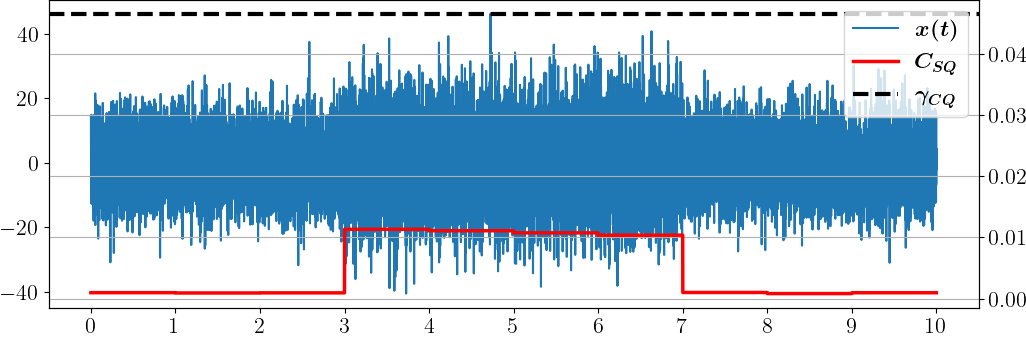}
    \vspace{0.1cm}
    \includegraphics[width = 12cm]{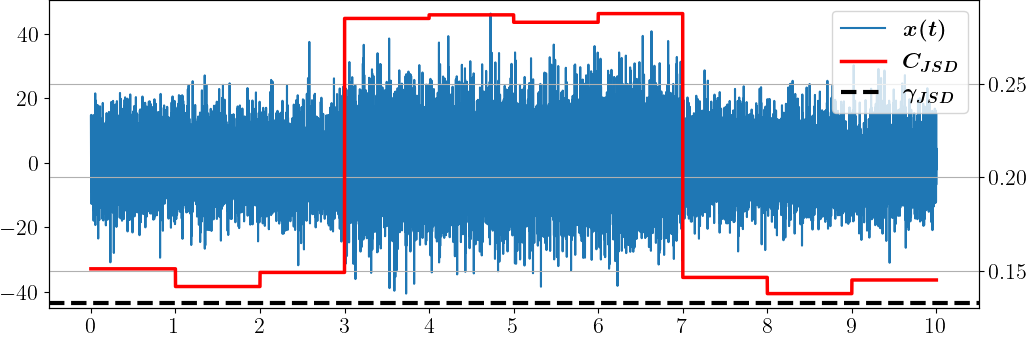}
    \vspace{0.1cm}
    \includegraphics[width = 12cm]{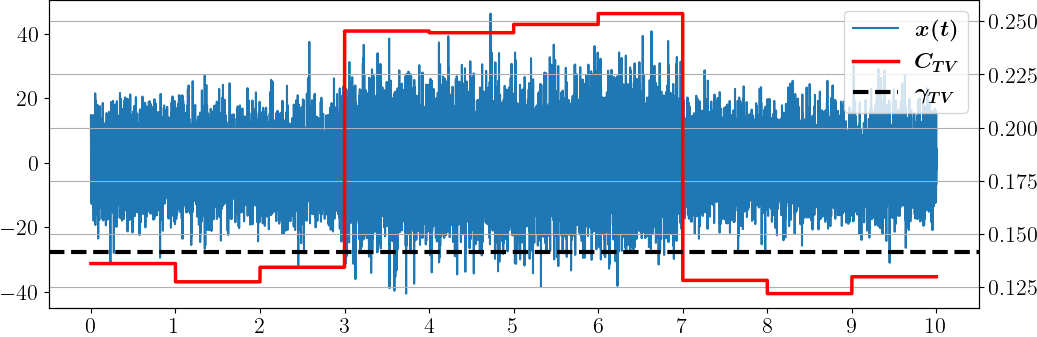}
    \caption{Thirty components, $K = 30$}
    \label{secondexp}
\end{figure}
\section{Conclusion}\label{Conclusion}

The article provides a theoretical justification for the usage of statistical complexity as a criterion for solving the problem of hypothesis testing when its error probability is close to one. Three variants of statistical complexity for different disequilibrium functions are considered. New notions of disequilibrium and statistical complexity based on the total variation measure are introduced. Information criteria are compared and the classes of discrete distributions which provide maximum for different types of statistical complexity are discovered. The values of maxima for fixed numbers of distribution samples are found. It is shown that the statistical complexity $C_{TV}$ based on total variation is directly related to the problem of hypothesis testing, while the statistical complexity $C_{JSD}$ based on Jensen-Shannon entropy gives a close estimate of $C_{TV}$ on sample distributions. In turn, $C_{SQ}$ is most promising for detecting an individual component over a uniform distribution. We propose a method for selecting the threshold for the decisive rule for the detection of a useful signal, taking into account the maximum values of the criteria obtained and show the effectiveness of this approach on the synthesized signals.

Future work will be devoted to the study of information criteria based on bivariate and multivariate distributions, as well as investigation of typical acoustic signals with realistic background noise.

\bmhead{Funding}

The work was supported by the Russian Science Foundation under grant \mbox{no 23-19-00134}.

\bibliography{sn-bibliography}



\end{document}